
%
%
%
%
%
\input amstex
\documentstyle{amsppt}
\loadbold
\def\cstar{$C^*$-algebra}
\magnification=\magstep 1

\topmatter
\title
Discretized CCR algebras
\endtitle

\author William Arveson
\endauthor

\affil Department of Mathematics\\
University of California\\Berkeley CA 94720, USA\\
29 July, 1991
\endaffil

\date {\it To the memory of John Bunce}
\enddate
\thanks
This research was supported in part by
NSF grant DMS89-12362
\endthanks
\abstract We discuss how the canonical commutation relations
must be modified in order to make appropriate numerical models
of quantum systems. The \cstar s
associated with the discretized $CCRs$ are the non-commutative
spheres of Bratteli, Elliott, Evans and Kishimoto.
\endabstract
\endtopmatter
\vfill\eject
\pagebreak
\document

\subheading{1. Introduction}

We consider the problem of discretizing the Hamiltonian of a
one-dimensional quantum system in a form that is appropriate
for carrying out numerical studies.  Specifically, we start with
a formal Schr\"odinger operator
$$
H = \frac{1}{2}P^2 + v(Q)
$$
acting on the Hilbert space $L^2(\Bbb R)$, where $P$ and $Q$
are the canonical operators
$$
\align
P &= -i\frac{d}{dx}\\
Q &= \text{multiplication by } x,
\endalign
$$
and $v$ is a real-valued continuous function of a real variable.
The problem of discretizing $H$ is that of finding an approximation
to $H$ which satisfies two requirements: $(a)$ the basic principles
of numerical analysis are satisfied, and $(b)$ the uncertainty
principle is preserved.

In \cite{1, \S\S\ 1--2}, we argued that in order to satisfy these
two conditions one must first replace $P$, $Q$ with the pair
$$
\align
P_\tau &= \frac{1}{\tau}\sin(\tau P)\\
Q_\tau &= \frac{1}{\tau}\sin(\tau Q).
\endalign
$$
Here, $\tau$ is a fixed positive real number, the numerical step size.
The discretized Hamiltonian is then defined as the following bounded
self-adjoint operator on $L^2(\Bbb R)$:
$$
H_\tau = \frac{1}{2}P_\tau^2 + v(Q_\tau).
$$

Obviously, $H_\tau$ belongs to the unital \cstar\ $C^*(P_\tau,Q_\tau)$
generated by $P_\tau$ and $Q_\tau$.  We show that when $\tau^2/\pi$
is irrational (e.g., when $\tau$ is a rational number), $C^*(P_\tau,Q_\tau)$
is isomorphic to the non-commutative sphere $\Cal B_{\tau^2}$ of Bratteli,
Evans, Elliott and Kishimoto \cite{5}\cite{6}; hence it is a simple \cstar\
with a unique trace.  We also describe the way in which the canonical
commutation relations must be ``discretized" in order to accommodate
pairs of operators $(P_\tau,Q_\tau)$ of this type. Together, these
observations serve to make a more philosophical point, namely
{\it non-commutative spheres will arise in any serious attempt to
model quantum systems on a computer}.

In the ``linear" case where $v$ has the form $v(x) = cx^2/2$, $c$
being a positive constant,
the operator $H_\tau$ turns out to be unitarily equivalent to an operator
of the form $\lambda M+\mu I$, where $\lambda$ and $\mu$ are real
constants and $M$ is the almost Mathieu Hamiltonian
$$
M = U + U^* + c(V + V^*),
$$
associated with a pair of unitary operators $U$, $V$ satisfying
$$
VU = e^{i4\tau^2}UV.
$$
An extensive amount of work has been done to compute the spectra
of such operators. Here, we mention only \cite{2}, \cite{3}, \cite{4},
\cite{9}, \cite{16} and refer the reader to the monograph \cite{8}
for further references.

Finally, I would like to thank Larry Schweitzer for pointing out
the references \cite{11} and \cite{13} (as well as the relevance of
his own work \cite{17}) in connection with the spectral invariance
property of the Banach $*$-algebra $\l^1(\Bbb Z\oplus\Bbb Z,\omega)$.

\subheading{2. Discretized CCR algebras}

Let $\theta$ be a real number such that $\theta/\pi$ is irrational,
and let $\omega$ be the bicharacter of the discrete abelian group
$G = \Bbb Z\oplus\Bbb Z$ defined by
$$
\omega((m,n),(p,q)) = e^{i(np-mq)\theta/2}\leqno{2.1}.
$$
A uniformly bounded family $\{D_x: x\in G\}$ of self-adjoint
operators on a Hilbert space $H$ is said to satisfy the
{\it discretized canonical commutation relations} if
$$
D_xD_y = \omega(x,y)D_{x+y}+\omega(y,x)D_{x-y},\qquad x,y\in G.\leqno{2.2}
$$

\remark{Remarks}Notice that $(2.2)$ is a generalization of
the elementary trigonometric identity
$$
2\cos A\cos B = \cos(A+B) + \cos(A-B),
$$
in which phase shifts have been added by way of the cocycle $\omega$.
Indeed, for any pair of real numbers $\alpha, \beta$, the function
$D: G\to\Bbb R$ defined by
$$
D(m,n) = 2\cos(\alpha m+\beta n)
$$
satisfies $(2.2)$ for the trivial cocycle $\omega = 1$.  It is related
to formula $(2.2)$ of \cite{5}, except that our operators are
self-adjoint and the phase factor is associated with a {\it nondegenerate}
bicharacter $\omega$.
\endremark

The purpose of this section is to associate a \cstar\ with the relations
$(2.2)$, and to point out some of its basic properties.  Let $\{D_x:x\in G\}$
satisfy $(2.2)$.  It is clear that the norm closed linear span
$$
\Cal D = \overline{span}\{D_x: x\in G\}
$$
is a separable \cstar.  Thus by passing from $H$ to the subspace $[\Cal DH]$
if necessary, we can assume that $\Cal D$ is nondegenerate.

\proclaim{Proposition 2.3}
\roster
\item"{(i)}" $D_0 = 2I.$
\item"{(ii)}" $D_{-x} = D_x.$
\item"{(iii)}" $\|D_x\|\leq 2$, for every $x\in G$.
\endroster
\endproclaim
\demo{proof}Setting $y=0$ in $(2.2)$ we obtain
$D_xD_0=2D_x$ for all $x\in G$, from which (i) is evident.
Setting $x=0$ in $(2.2)$ now leads to
$2D_y = D_0D_y=D_y+D_{-y}$, hence (ii).  For (iii), let
$$
M = \sup_{x\in G}\|D_x\|.
$$
By hypothesis, $M<\infty$.  Moreover, setting $y=x$ in $(2.2)$ gives
$$
D_x^2 = D_{2x}+D_0 = D_{2x} + 2I
$$
and thus $M^2 \leq M+2$.  This inequality implies that
$-1\leq M\leq 2$, hence (iii)\qed\enddemo

We now construct a Banach $*$-algebra whose representations are
associated with operator realizations of $(2.2)$.  Let $\l^1(G,\omega)$
denote the Banach space of all absolutely summable complex functions on $G$,
endowed with the multiplication and involution
$$
\align
f*g(x) &= \sum_{y}\omega(y,x)f(y)g(x-y)\\
f^*(x) &= \overline{f(-x)}.
\endalign
$$
It is easily checked that the linear subspace
$$
D_\theta = \{f\in \l^1(G,\omega): f(-x) = f(x),\quad x\in G\}
$$
is in fact a $*$-subalgebra of $\l^1(G,\omega)$.  Of course, the adjoint
operation in $D_\theta$ simplifies to $f^*(x) = \overline{f(x)}$.  Moreover,
$D_\theta$ is linearly spanned by the elements
$$
d_x = \delta_x + \delta_{-x},
$$
$\delta_x$ denoting the unit function supported at $x$, and one has
$$
\align
d_xd_y &=\omega(x,y)d_{x+y} + \omega(y,x)d_{x-y}\\
\|d_x\| &= 2\\
d_x &=d_{-x} = d_x^*.
\endalign
$$

\proclaim{Proposition 2.4} Let $\{D_x: x\in G\}$ be a uniformly
bounded family of self-adjoint operators on a Hilbert space $H$
satisfying $(2.2)$.  Then there is a unique representation
$\pi : D_\theta\to \Cal B(H)$ such that
$$
\pi(d_x) = D_x,\qquad x\in G.
$$
\endproclaim
\demo{proof} By proposition $(2.3)$, we know that $\|D_x\| \leq 2$; hence
$$
\pi(f) = \frac{1}{2} \sum_{x\in G}f(x)D_x
$$
defines a contractive self-adjoint linear mapping of $D_\theta$ into
$\Cal B(H)$.  Moreover, using $(2.2)$ we have
$$
\align
\pi(f)\pi(g)&=\frac{1}{4}\sum_{x,y}f(x)g(y)(\omega(x,y)
D_{x+y} + \omega(y,x)D_{x-y})\\
&=\frac{1}{4}\sum_{z,x} f(x)g(z-x)\omega(x,z)D_z
+ \frac{1}{4}\sum_{z,x} f(x)g(x-z)\omega(-z,x)D_z.
\endalign
$$
Using the fact that $g(x-z)=g(z-x)$ and $\omega(-z,x)=\omega(x,z)$,
the right side becomes
$$
\frac{1}{2}\sum_z(\sum_x f(x)g(z-x)\omega(x,z))D_z = \pi(f*g),
$$
as required.

Finally, taking $f = \delta_x + \delta_{-x} = d_x$ and using $(ii)$ of $(2.3)$,
we find that
$$
\pi(d_x) = \frac{1}{2}(f(x)D_x + f(-x)D_x) = D_x
$$
as required \qed\enddemo

\remark{Remarks} It follows that the enveloping \cstar\ $C^*(D_\theta)$
is the universal \cstar\ generated by the commutation relations $(2.2)$.

Let $\alpha$ be an automorphism of the discrete abelian group
$\Bbb Z\oplus\Bbb Z$. Then $\alpha$ is given by a $2\times 2$
integer matrix
$$
\left(\matrix
a & b\\
c & d
\endmatrix\right)
$$
by way of $\alpha(m,n) = (am + bn, cm + dn)$,
where $a,b,c,d \in \Bbb Z$ satisfy the condition
$$
\det\alpha = ad - bc = \pm 1.
$$
It follows that
$$
\omega(\alpha x,\alpha y) =
\cases\omega(x,y),&\text{if $\det\alpha=+1$}\\
		\omega(y,x),&\text{if $\det\alpha=-1$}.\endcases
$$
Hence the group $SL(2,\Bbb Z)$ of determinant $1$ automorphisms acts
naturally on $D_\theta$ (resp. $C^*(D_\theta)$) as a group of
$*$-automorphisms.  Any $\alpha\in aut(\Bbb Z\oplus\Bbb Z)$ satisfying
$\det\alpha = -1$
gives rise to a $*$-anti-automorphism of $D_\theta$ (resp. $C^*(D_\theta)$).

Finally, notice that there is a natural $*$-homomorphism which carries
$D_\theta$ into the irrational rotation \cstar\ $\Cal A_\theta$.  Indeed,
$D_\theta$ is obviously contained in the larger Banach $*$-algebra
$\l^1(\Bbb Z\oplus\Bbb Z,\omega)$ obtained by simply dropping the
requirement that $f(-x) = f(x)$.  It is clear that
$\l^1(\Bbb Z\oplus\Bbb Z,\omega)$ is the universal Banach $*$-algebra
generated by unitary operators $\{W_x: x\in \Bbb Z\oplus\Bbb Z\}$ satisfying
$$
W_xW_y = \omega(x,y)W_{x+y},\qquad x,y\in\Bbb Z\oplus\Bbb Z.
$$
Because of the formula $(2.1)$ giving $\omega$ in terms of $\theta$,
the unitary elements $U, V$ defined by $U = W_{(1,0)}, V=W_{(0,1)}$
satisfy $VU = e^{i\theta} UV$, and of course they generate
$\l^1(\Bbb Z\oplus\Bbb Z,\omega)$ as a Banach $*$-algebra.  It follows
that the enveloping \cstar\ of $\l^1(\Bbb Z\oplus\Bbb Z,\omega)$ is
$\Cal A_\theta$.  Thus
we obtain a morphism of $D_\theta$ into $\Cal A_\theta$ by simply
restricting the completion map
$$
\gamma : \l^1(\Bbb Z\oplus\Bbb Z,\omega) \to \Cal A_\theta
$$
to $D_\theta$.  By the universal property of enveloping
\cstar s there is correspondingly a unique morphism of \cstar s
$$
\gamma_B : C^*(D_\theta) \to \Cal A_\theta.
$$
In the next section it will be shown that $\gamma_B$ is injective
and we will identify its range.
\endremark

\subheading{3. Spectral Invariance and Extensions of States}

Let $A$ be a Banach $*$-algebra with unit, and let
$$
A^+ = \overline{\{a_1^*a_1 + a_2^*a_2+\dots+a_n^*a_n :
a_k\in A, n \geq 1\}}
$$
denote the closed positive cone in $A$.  For simplicity, we
assume throughout this section that the completion map $\gamma$
of $A$ into its enveloping \cstar\ is {\it injective}.

Let $B$ be a unital self-adjoint Banach subalgebra of $A$.  We
are interested in determining  whether or not the \cstar\
obtained by closing $\gamma(B)$ in the norm of $C^*(A)$ is the
enveloping \cstar\ of $B$.  More precisely, we seek conditions
under which the $*$-homomorphism $\gamma_B : C^*(B) \to C^*(A)$
defined by the commutative diagram
$$
\CD
B	@>\text{incl}>>	A\\
@VVV													@VVV \\
C^*(B)	@>\gamma_B>>	C^*(A)
\endCD \tag{$3.1$}
$$
should be {\it injective}.  Elementary considerations show that
the following three conditions are equivalent:

\roster
\item $\gamma_B$ is injective.
\item Every positive linear functional on $B$ can be extended
to a positive linear functional on $A$.
\item $A^+\cap B \subseteq B^+$.
\endroster
Note, for example, that the implication $(3) \implies (2)$ is
the extension theorem of M. G. Krein \cite{15, p. 227},
whereas $(2) \implies (3)$ follows from a standard separation theorem.
It is not hard to find examples showing that these conditions
are not always satisfied (see Appendix).

$A$ is said to have the {\it spectral invariance property} if for
every element $a\in A$ which is invertible in $C^*(A)$, we have
$a^{-1}\in A$.  This is equivalent to the assertion that the
spectrum of any element of $A$ is the same whether it is computed
in $A$ or in $C^*(A)$, or that $A$ is closed under the holomorphic
functional calculus of $C^*(A)$ (see \cite{10, p. 52} for
further significant consequences of spectral invariance in more
general Fr\`echet algebras).

A familiar Tauberian theorem of Wiener asserts that if a continuous
function on the unit circle never
vanishes and has an absolutely convergent Fourier series
$$
f(e^{i\theta}) = \sum_{n = -\infty}^{+\infty}a_ne^{in\theta},
$$
$\sum|a_n| < \infty$, then $1/f$
has an absolutely convergent Fourier series.  Of course, this
is precisely the assertion that the group
algebra $\l^1(\Bbb Z)$ has the spectral invariance property. While
this theorem has a simple proof using the Gelfand theory, it is
certainly not a triviality.

The significance of spectral invariance for our purposes derives
from the following.

\proclaim{Proposition 3.2}Let $A$ be a unital Banach $*$-algebra which
admits spectral invariance.  Then for every self-adjoint unital Banach
subalgebra $B$ of $A$, the natural $*$-homomorphism
$$
\theta_B: C^*(B) \to C^*(A)
$$
is injective.
\endproclaim
\demo{proof}We will verify property $(3)$ above by showing that
$A^+\cap B \subseteq B^+$.  We may clearly assume that
$A\subseteq C^*(A)$, as a self-adjoint subalgebra which is a
Banach algebra relative to a larger norm than that of $C^*(A)$.

Choose $x\in A^+\cap B$; without loss of generality we may
assume that the $B$-norm of $x$ is less than $1$.  Since $x$ belongs
to the positive cone of $C^*(A)$ its spectrum in $C^*(A)$ is nonnegative.
By spectral invariance we have $\sigma_A(x) \subseteq [0,1)$.  Moreover,
since $\sigma_A(x)$ cannot separate the complex plane, we see from
the spectral permanence theorem that
$\sigma_B(x) = \sigma_A(x)\subseteq [0,1)$.  Hence for sufficiently
small $\epsilon$ we have $\sigma_B(x + \epsilon1) \subseteq (\epsilon,1)$.
Thus we may apply the power series
$$
\sqrt t = \sum_{n = 0}^\infty a_n(1-t)^n,\qquad |1 - t| < 1
$$
to the element $x + \epsilon1$ to obtain a square root in B, i.e., a
self-adjoint element $h\in B$ satisfying $x + \epsilon1 = h^2$.
This shows
that $x + \epsilon1 \in B^+$, and we obtain the desired  conclusion by
allowing $\epsilon$ to tend to zero \qed
\enddemo

We now apply this to show that the enveloping \cstar\ $C^*(D_\theta)$
is isomorphic to the non-commutative sphere $\Cal B_\theta$ of \cite6.
If we realize the irrational rotation \cstar\ $\Cal A_\theta$ as the
\cstar\ generated by a pair of unitary operators $U$, $V$ satisfying
$VU = e^{i\theta}UV$, then there is a unique automorphism $\sigma$ of
$\Cal A_\theta$ satisfying $\sigma(U) = U^{-1}$, $\sigma(V) = V^{-1}$.
In case $\theta/\pi$ is irrational, $\Cal B_\theta$ is defined to be the
fixed subalgebra
$$
\Cal B_\theta = \{a\in\Cal A_\theta : \sigma(a) = a\}.
$$
Let $\{W_x: x\in \Bbb Z\oplus\Bbb Z\}$ be the family of unitary
operators in $\Cal A_\theta$ defined by
$$
W_{(m,n)}=e^{imn\theta/2}U^mV^n,\qquad m,n\in \Bbb Z.
$$
One verifies easily that
$$
W_xW_y = \omega(x,y)W_{x+y},
$$
where $\omega$ is the bicharacter on $\Bbb Z\oplus\Bbb Z$ defined
in (2.1), and moreover the action of $\sigma$ is given by
$$
\sigma(W_x) = W_{-x},\qquad x\in \Bbb Z\oplus\Bbb Z.
$$
Since $\Cal A_\theta$ is spanned by $\{W_x: x\in \Bbb Z\oplus\Bbb Z\}$,
we conclude that $\Cal B_\theta$ is spanned by
$\{W_x+W_{-x}: x\in \Bbb Z\oplus\Bbb Z\}$.

\proclaim{Corollary}Suppose $\theta$ is not a rational multiple of
$\pi$, and let $\alpha: D_\theta\to\Cal A_\theta$
be the morphism defined by
$$
\alpha(d_x) = W_x + W_{-x},\qquad x\in\Bbb Z\times\Bbb Z.
$$
Then the natural extension $\tilde\alpha: C^*(D_\theta)\to \Cal A_\theta$
gives an isomorphism of \cstar s
$$
C^*(D_\theta) \cong\Cal B_\theta.
$$
\endproclaim
\demo{proof}
Let $G = \Bbb Z\oplus\Bbb Z$ and let $\omega:G\times G\to\Bbb T$
be the bicharacter of (2.1).  Consider the Banach
$*$-algebra $\l^1(G, \omega)$, where multiplication and involution
are defined respectively by
$$
\align
f*g(x) &= \sum_{y}\omega(y,x)f(y)g(y-x)\\
f^*(x) &= \overline{f(-x)}.
\endalign
$$

Notice first that $C^*(G,\omega)$ is naturally identified with the
irrational rotation \cstar\ $A_\theta = C^*(U,V)$, where $U$ and $V$
are unitary operators satisfying the above relation
$VU = e^{i\theta}UV$.  Indeed, letting $\{W_x: x\in G\}$ be the
operators of $\Cal A_\theta$ defined in the preceding remarks,
it is clear that we can define a morphism of $\l^1(G,\omega)$ into
$\Cal A_\theta$ by
$$
\gamma(\delta_x) = W_x,\qquad x\in G,
$$
$\delta_x$ denoting the unit function at $x$.  The range of $\gamma$ is
dense in $\Cal A_\theta$, and the natural extension of $\gamma$ to
$C^*(\l^1(G,\omega))$ is injective because of the familiar universal
property of such pairs $U$, $V$.

$D_\theta$ is indentified (via an isometric isomorphism of
Banach $*$-algebras) with a
subalgebra of  $\l^1(G,\omega)$,
$$
D_\theta = \{f\in\l^1(G,\omega) : \sigma_0(f) = f\},
$$
where $\sigma_0$ is the $*$-automorphism of $\l^1(G,\omega)$ given by
$$
\sigma_0(f)(x) = f(-x),\qquad x\in G.
$$
It is clear that the restriction of $\gamma$ to $D_\theta$ carries
$d_x = \delta_x + \delta_{-x}$ to
$W_x + W_{-x}$; hence by the preceding remarks
$\gamma(D_\theta)$ is a dense $*$-subalgebra
of $\Cal B_\theta$.  Thus $\gamma$ extends naturally to a surjective
$*$-homomorphism of $C^*(\l^1(G,\omega))$ onto $\Cal B_\theta$, and
it remains only to show that the latter morphism is injective.

Now it is known that $\l^1(G,\omega)$ admits spectral invariance
(see \cite{13, Satz 5} for example, or apply Theorem 1.1.3 of \cite{17}
together with the results of \cite{11} on the symmetry of the group algebra
of the rank 3 discrete Heisenberg group); hence Proposition 3.2
implies that $\gamma\restriction_{D_\theta}$ extends uniquely to a
$*$-isomorphism of
$C^*(D_\theta)$ onto $\overline{\gamma(D_\theta)} = \Cal
D_\theta$ \qed \enddemo

\subheading{4. Representations}

In this section we make some general comments about the representation
theory of the discretized CCRs $(2.2)$.  We assume throughout that
$\theta$ is a real number such that $\theta/\pi$ is irrational.

\remark{Remark 4.1: Finite representations}

The unique trace on the irrational rotation
algebra $\Cal A_\theta$ gives rise to a representation of $\Cal A_\theta$
which generates the hyperfinite $II_1$ factor $R$.  The closure of
$\Cal B_\theta$ in this representation is a sub von Neumann algebra
of $R$.  Since $\Cal B_\theta$ has a unique tracial state \cite{5},
it follows that the closure of $\Cal B_\theta$ is a subfactor of $R$,
and hence is also isomorphic to $R$.  Moreover, since
$\Cal B_\theta$ is also simple \cite{5}, any
finite representation of $\Cal B_\theta$ is quasi-equivalent to
this one.

It is not hard to show that the subfactor of $R$ generated by $\Cal B_\theta$
in the above representation has Jones index 2.  Since any two subfactors
of $R$ of index 2 are known to be isomorphic \cite{12}, we have here
a very stable invariant for the embedding of the discretized $CCR$
algebra in the irrational rotation algebra $\Cal A_\theta$.

In particular, by the corollary of 3.2 we may conclude from these
remarks that {\it there
is a representation of the discretized $CCR$s $(2.2)$ which generates
$R$ as a von Neumann algebra; moreover any finite representation
of the discretized $CCR$s is quasi-equivalent to this one}.
\endremark

Now let $\tau$ be a positive real number such that $\tau^2/\pi$ is irrational,
and let $P_\tau$ and $Q_\tau$ be the discretized canonical operators on
$L^2(\Bbb R)$ associated with the step size $\tau$ as in section 1.
We want to make explicit the relation that exists between the pair
$(P_\tau,Q_\tau)$ and the \cstar\ $C^*(\Cal B_{\tau^2})$ discussed in
section 2.

\proclaim{Theorem 4.2}There is a unique representation $\pi$
of $D_{\tau^2}$ on $L^2(\Bbb R)$ satisfying
$$
\align
\pi(d_{(1,0)}) &= 2\tau Q_\tau,\\
\pi(d_{(0,1)}) &= 2\tau P_\tau.
\endalign
$$
$\pi(D_{\tau^2})$ and $\{P_\tau, Q_\tau\}$ generate the same unital
\cstar.  Thus, the three \cstar s
$$
C^*(D_{\tau^2}), \quad\Cal B_{\tau^2},\quad C^*(P_\tau,Q_\tau)
$$
are mutually isomorphic.
\endproclaim

\demo{proof} Let $U$, $V$ be the one-parameter groups
$$
\align
U_tf(x) &= e^{itx}f(x),\\
V_tf(x) &= f(x+t)\qquad f\in L^2(\Bbb R).
\endalign
$$
As in section 1 we have
$$
\aligned
Q_\tau &= \frac{1}{2i\tau}(U_\tau - U_{-\tau})=\frac{1}{\tau}\sin(\tau Q),\\
P_\tau &= \frac{1}{2i\tau}(V_\tau - V_{-\tau})=\frac{1}{\tau}\sin(\tau P).
\endaligned\leqno{4.3}
$$

We claim first that the {\it sines} in $(4.3)$ can be replaced by {\it cosines}
in the sense that the pair $(P_\tau,Q_\tau)$ is unitarily equivalent to the
pair $(\tilde P_\tau,\tilde Q_\tau)$ given by
$$
\aligned
\tilde Q_\tau &= \frac{1}{2\tau}(U_\tau + U_{-\tau})\\
\tilde P_\tau &= \frac{1}{2\tau}(V_\tau + V_{-\tau}).
\endaligned\leqno{4.4}
$$
To see this, put $\lambda = \pi/{2\tau}$ and let $R$ denote the reflection
on $L^2(\Bbb R)$ given by $Rf(x)=f(-x)$.  Consider the unitary operator
$$
W = RU_{-\lambda}V_\lambda.
$$
Using the commutation relations $V_tU_s = e^{ist}U_sV_t$ together with
$RU_sR^* = U_{-s}$ and $RV_tR^* = V_{-t}$, one finds that
$$
\align
WU_sW^* &= e^{i\lambda s}U_{-s}\\
WV_tW^* &= e^{i\lambda t}V_{-t}.
\endalign
$$
Noting that $e^{i\lambda\tau}=\sqrt{-1}$, we obtain $(4.4)$ by applying
$ad W$ to $(4.3)$, i.e.,
$$
\align
W Q_\tau W^* &= \frac{1}{2\tau}(U_\tau + U_{-\tau})=\frac{1}{2\tau}\cos(\tau
Q)\\
W P_\tau W^* &= \frac{1}{2\tau}(V_\tau + V_{-\tau})=\frac{1}{2\tau}\cos(\tau
P).
\endalign
$$
We may therefore assume that the pair $(Q_\tau,P_\tau)$ is defined by $(4.4)$.

For each $x = (m,n)\in \Bbb Z\oplus\Bbb Z$, define a unitary operator $W_x$ by
$$
W_{(m,n)}=e^{imn\tau^2/2}U_{m\tau}V_{n\tau}.
$$
A straightforward computation shows that the family of unitaries
$\{W_x: x\in\Bbb Z\oplus\Bbb Z\}$ satisfies
$$
W_xW_y = \omega(x,y)W_{x+y}
$$
$\omega$ being the cocycle of $(2.1)$ for the value $\theta=\tau^2$, and
hence there is a representation $\pi$ of $\l^1(\Bbb Z\oplus\Bbb Z,\tau^2)$
on $L^2(\Bbb R)$ such that
$$
\pi(w_x) = W_x,\qquad x\in\Bbb Z\oplus\Bbb Z.
$$
It is clear that $\pi$ carries $d_{(1,0)}$ (resp. $d_{(0,1)}$) to
$U_\tau + U_{-\tau}=2\tau Q_\tau$ (resp. $2\tau P_\tau$).

It remains to show that the restriction of $\pi$ to
$C^*(D_{\tau^2})$ is uniquely defined by its values on the
two elements $d_{(1,0)}, d_{(0,1)}\}$, and that $Q_\tau$ and $P_\tau$ generate
$\pi(C^*(D_{\tau^2}))$ as a unital \cstar.  We will prove both
by showing that the two elements $\{d_{(1,0)}, d_{(0,1)}\}$ and the
identity generate
the Banach $*$-algebra $D_{\tau^2}$.  It is not hard to adapt the results
of \cite{5} to prove that these three elements generate $D_{\tau^2}$.
Instead, we present the following argument since it gives somewhat more
structural information.

Actually, we will give a fairly explicit method for calculating each element
$d_x=\delta_x+\delta_{-x}$ in terms of the self-adjoint elements
$p=d_{(1,0)}$ and $q=d_{(0,1)}$, using a
``generating function" for the family $\{d_x: x\in \Bbb Z\oplus\Bbb Z\}$.
Indeed, it suffices to establish the following lemma.

\proclaim{Lemma 4.5} Let $\theta$ be a real number such that $\theta/\pi$
is irrational, and consider the real-analytic function
$F: (-1,1)\times(-1,1)\to D_\theta$ defined by
$$
F(s,t) = \sum_{m,n=-\infty}^{+\infty}
											s^{|m|}t^{|n|}e^{-imn\theta/2}d_{(m,n)}.\leqno{4.6}
$$
\roster
\item"{(i)}"
For $-1<u<1, -2\leq x\leq 2$, let
$$
\phi(u,x) =\frac{1-u^2}{1+u^2 -ux}.
$$
Noting that $\phi$ is separately analytic in each variable, we have
$$
F(s,t) = 2\phi(s,q)\phi(t,p),\qquad |s|,|t| <1,
$$
where $q,p$ are the elements of $D_\theta$ defined by
$$
q=d_{(1,0)}, \qquad p=d_{(0,1)}.
$$
\item"{(ii)}"
The Banach $*$-algebra $D_{\theta}$ is spanned by the set $\Cal F\cup\Cal F^*$,
where
$$
\Cal F = \{F(s,t): |s|,|t| < 1\}.
$$
\endroster\endproclaim

{proof of $(i)$}  Let $\omega$ be the bicharacter of $\Bbb Z\oplus\Bbb Z$
defined by
$$
\omega((p,q),(m,n))=e^{i(qm-pn)\theta/2}
$$
and let $u,v$ be the following elements of $\l^1(\Bbb Z\oplus\Bbb Z,\omega)$:
$$
u=\delta_{(1,0)},\qquad v=\delta_{(0,1)}.
$$
Then $w_{(m,n)}=e^{imn\theta/2}u^mv^n$, hence
$$
d_{(m,n)}= e^{imn\theta/2}(u^mv^n+u^{-m}v^{-n}).
$$
It follows that
$$
\align
F(s,t) &= \sum_{m,n={-\infty}}^{\infty}s^{|m|}t^{|n|}(u^mv^n+u^{-m}v^{-n})
					=2\sum_{m,n={-\infty}}^{\infty}s^{|m|}t^{|n|}u^mv^n\\

&=2\sum_{m={-\infty}}^{\infty}s^{|m|}u^m\sum_{n={-\infty}}^{\infty}t^{|n|}v^n.
\endalign
$$
An elementary calculation shows that if $z$ is any complex number having
absolute value $1$ and $-1<s<1$, then
$$
\sum_{m={-\infty}}^{\infty}s^{|m|}z^m =
\frac{1-s^2}{1+s^2 -s(z+\bar z)}=\phi(s,z+\bar z).
$$
Since $q = d_{(1,0)}=u+u^*$ and $p = d_{(0,1)} = v+v^*$, the assertion $(i)$
follows
from the analytic functional calculus.

To prove $(ii)$, let $A_{pq}$ be the coefficients in the power series
expansion of $F$,
$$
F(s,t)=\sum_{p,q=0}^{\infty}A_{pq}s^pt^q.
$$
Obviously, $\{F(s,t): s,t\in (-1,1)\}$ and $\{A_{pq}:p,q\geq0\}$ have
the same closed linear span.  Using the fact that $d_{(-m,-n)}=d_{(m,n)}$, a
straightforward computation shows that
$$
A_{pq} = 2e^{-ipq\theta/2}d_{(p,q)} + 2e^{ipq\theta/2}d_{(-p,q)}.
$$
Thus,
$$
\align
A_{0q} &= 2(d_{(0,q)} + d_{(0,q)})=4d_{(0,q)},\qquad {\text and}\\
A_{p0} &= 2(d_{(p,0)} + d_{(-p,0)})=4d_{(p,0)}.
\endalign
$$
In the remaining cases where $pq\neq 0$, the determinant of the coefficients
of the $2\times 2$ system of operator equations
$$
\align
A_{pq} &= 2e^{-ipq\theta/2}d_{(p,q)}+ 2e^{ipq\theta/2}d_{(-p,q)})\\
A_{pq}^* &= 2e^{ipq\theta/2}d_{(p,q)}+ 2e^{-ipq\theta/2}d_{(-p,q)})\tag{4.7}
\endalign
$$
is $4(e^{-ipq\theta}-e^{ipq\theta})\neq 0$, and in particular we can
solve $(4.7)$ for $d_{(p,q)}$ as a complex linear combination of
$A_{pq}$ and $A_{pq}^*$.  This argument shows that the closed linear span
of $\Cal F\cup\Cal F^*$ contains $\{d_{(p,q)}: p,q\in \Bbb Z\}$, and $(ii)$
follows.  That completes the proof of Theorem 4.2 \qed
\enddemo

\remark{Remark} In some very recent work \cite{7}, Bratteli and Kishimoto
have established the striking result that $\Cal B_\theta$ is an $AF$-algebra.
\endremark
\newpage

\subheading{Appendix: Failure of Extensions}We present a
simple example of a pair of commutative unital Banach $*$-algebras
$B\subseteq A$ such that $A$ is a subalgebra of its
enveloping \cstar, but such that the natural morphism
$\gamma_B:C^*(B) \to C^*(A)$
is not injective.  Let $A$ be the
algebra of all complex-valued continuous functions defined
on the annulus $\{1\leq|z|\leq 2\}$ which are analytic in its
interior.  With norm and involution defined by
$$
\|f\| = \sup_{1\leq|z|\leq 2}|f(z)|,\qquad
f^*(z) = \bar f(\bar z),
$$
$\bar f$ denoting the complex conjute of $f$, A is a unital Banach
$*$-algebra.  $C^*(A)$ is the commutative \cstar\ $C(X)$,
$$
X=[-2,-1]\cup[+1,+2]
$$
denoting the intersection of the annulus $\{1\leq|z|\leq 2\}$
with the real axis, and the completion map
$\gamma : A\to C(X)$ is defined by restriction to $X$.  Let
$B$ be the norm closure of all holomorphic polynomials in $A$.  Then
$B$ is a self-adjoint subalgebra whose enveloping \cstar\ is $C(Y)$,
$Y$ being the intersection of the polynomially convex hull of the
annulus with the real axis, namely
$$
Y = [-2,+2].
$$

The morphism $\gamma_B:C(Y) \to C(X)$ is given by restriction to
$X$, and hence there is a nontrivial kernel.  Put differently, for
every real $\lambda\in (-1,+1)$, the complex homomorphism of $B$ defined
by
$$
\omega_\lambda(f) = f(\lambda), \qquad f\in B
$$
is a bounded positive linear functional on $B$ which cannot be extended
to a positive linear functional on $A$.

\end